\documentclass{amsart}
\usepackage{amssymb,latexsym}
\newcommand{\pq}[1]{[\,#1\,]_{p,q}}

\numberwithin{equation}{section}
\newtheorem{theorem}{Theorem}[section]

\newtheorem{corollary}[theorem]{Corollary}

\newtheorem{remark}[theorem]{Remark}
\newtheorem{lemma}[theorem]{Lemma}

\newtheorem{example}[theorem]{Example}

\begin{document}

\pagenumbering{arabic} \pagestyle{headings}
\def\sof{\hfill\rule{2mm}{2mm}}
\def\ls{\leq}
\def\gs{\geq}
\def\SS{\mathcal S}
\def\qq{{\bold q}}
\def\txx{{\frac1{2\sqrt{x}}}}
\def\mn{\mbox{-}}

\title{Restricted permutations by patterns of type $(2,1)$}

\author{Toufik Mansour}
\maketitle

\begin{center}{
LaBRI (UMR 5800), Universit\'e Bordeaux 1, 351 cours de la Lib\'eration,\\
33405 Talence Cedex, France\\[4pt]
{\tt toufik@labri.fr}
}\end{center}

\section*{Abstract}
Recently, Babson and Steingrimsson (see \cite{BS}) introduced
generalized permutations patterns that allow the requirement that
two adjacent letters in a pattern must be adjacent in the
permutation.

In this paper we study the generating functions for the number of
permutations on $n$ letters avoiding a generalized pattern $ab\mn
c$ where $(a,b,c)\in S_3$, and containing a prescribed number of
occurrences of generalized pattern $cd\mn e$ where $(c,d,e)\in
S_3$. As a consequence, we derive all the previously known results
for this kind of problems, as well as many new results.
\section{Introduction}
{\bf Classical patterns.} Let $\alpha\in S_n$ and $\tau\in S_k$ be
two permutations. We say that $\alpha$ {\em contains} $\tau$ if
there exists a subsequence $1\leq i_1<i_2<\cdots<i_k\leq n$ such
that $(\alpha_{i_1},\dots,\alpha_{i_k})$ is order-isomorphic to
$\tau$; in such a context $\tau$ is usually called a {\em
pattern}. We say $\alpha$ {\em avoids} $\tau$, or is $\tau$-{\em
avoiding}, if such a subsequence does not exist. The set of all
$\tau$-avoiding permutations in $S_n$ is denoted $S_n(\tau)$. For
an arbitrary finite collection of patterns $T$, we say that
$\alpha$ avoids $T$ if $\alpha$ avoids any $\tau\in T$; the
corresponding subset of $S_n$ is denoted $S_n(T)$.

While the case of permutations avoiding a single pattern has
attracted much attention, the case of multiple pattern avoidance
remains less investigated. In particular, it is natural, as the
next step, to consider permutations avoiding pairs of patterns
$\tau_1$, $\tau_2$. This problem was solved completely for
$\tau_1,\tau_2\in S_3$ (see~\cite{SS}), for $\tau_1\in S_3$ and
$\tau_2\in S_4$ (see~\cite{W}), and for $\tau_1,\tau_2\in S_4$
(see~\cite{B1,Km} and references therein). Several recent papers
deal with the case $\tau_1\in S_3$, $\tau_2\in S_k$ for various
pairs $\tau_1,\tau_2$ (see~\cite{CW,Kr,MV3} and references
therein). Another natural question is to study permutations
avoiding $\tau_1$ and containing $\tau_2$ exactly $t$ times. Such
a problem for certain $\tau_1,\tau_2\in S_3$ and $t=1$ was
investigated in~\cite{R}, and for certain $\tau_1\in S_3$,
$\tau_2\in S_k$ in~\cite{RWZ,MV1,Kr,MV3,MV2,MV4}. The tools
involved in these papers include continued fractions, Chebyshev
polynomials, and Dyck paths.\\

{\bf Generalized patterns.} In~\cite{BS} Babson and Steingrimsson
introduced generalized permutation patterns that allow the
requirement that two adjacent letters in a pattern must be
adjacent in the permutation. The idea for Babson and Steingrimsson
in introducing these patterns was study of Mahonian statistics.

Following~\cite{C}, we define our {\em generalized patterns} as
words on the letters $1,2,3,\dots$ where two adjacent letters may
or may not be separated by a dash. The absence of a dash between
two adjacent letters in a pattern indicates that the corresponding
letters in the permutation must be adjacent, and in the order
(order-isomorphic) given by the pattern. For example, the subword
$23\mn1$ of a permutation $\pi=(\pi_1,\pi_2,\cdots,\pi_n)$ is a
subword $(\pi_i,\pi_{i+1},\pi_{j})$ where $i+1<j$ such that
$\pi_j<\pi_{i}<\pi_{i+1}$. We say $\tau$ a generalized pattern of
type $(2,1)$ if it has the form $ab\mn c$ where $(a,b,c)\in S_3$.

\begin{remark}
There exist six generalized patterns of type $(2,1)$ which are
$12\mn3$, $13\mn2$, $21\mn3$, $23\mn1$, $31\mn2$, and $32\mn1$. By
the complement symmetric operation we get three different classes:
$\{12\mn3, 32\mn1\}$, $\{13\mn2,31\mn2\}$, and $\{21\mn3,
23\mn1\}$.
\end{remark}

While the case of permutations avoiding a single pattern has
attracted much attention, the case of multiple pattern avoidance
remains less investigated. In particular, it is natural, as the
next step, to consider permutations avoiding pairs of generalized
patterns $\tau_1$, $\tau_2$. This problem was solved completely
for $\tau_1,\tau_2$ two generalized patterns of length three with
exactly one adjacent pair of letters (see \cite{CM1}). Claesson
and Mansour \cite{CM2} showed (by Clarke, Steingr{\'\i}msson and
Zeng \cite[Corollary~11]{ClStZe97} result) the distribution of the
patterns $2\mn31$ and $31\mn2$ is given by Stieltjes continued
fraction as following:

\begin{theorem}\label{CSZ}
  The following Stieltjes continued fraction expansion holds
  $$
  \sum\limits_{n\geq 0}\sum\limits_{\pi\in S_n} x^{1+(12)\pi } y^{(21)\pi}
  p^{(2\mn31)\pi} q^{(31\mn2)\pi} t^{|\pi|} =
  \cfrac{1}{1 -
    \cfrac{x \pq 1 t}{1 -
      \cfrac{y \pq 1 t}{1 -
        \cfrac{x \pq 2 t}{1 -
          \cfrac{y \pq 2 t}{\quad\ddots}}}}}
  $$
  where $\pq n = q^{n-1}+pq^{n-2}+\cdots+p^{n-2}q+p^{n-1}$,
  $(\tau)\pi$ is the number of occurrences of $\tau$ in $\pi$.
\end{theorem}

In the present paper, as consequence to \cite{CM1} (see also
\cite{C,CM2}), we exhibit a general approach to study the number
of permutations avoiding a generalized pattern of type $(2,1)$ and
containing a prescribed number of occurrences of $\tau$ a
generalized pattern of type $(2,1)$. As a consequence, we derive
all the previously known results for this kind of problems, as
well as many new results.
\section{Avoiding $12\mn3$}
Let $f_{\tau;r}(n)$ be the number of all permutations in
$S_n(12\mn3)$ containing $\tau$ exactly $r$ times. The
corresponding exponential and ordinary generating function we
denote by $\mathcal{F}_{\tau;r}(x)$ and $F_{\tau;r}(x)$
respectively; that is,
$\mathcal{F}_{\tau;r}(x)=\sum\limits_{n\geq0}\frac{f_{\tau;r}(n)}{n!}x^n$
and $F_{\tau;r}(x)=\sum\limits_{n\geq0}f_{\tau;r}(n)x^n$. The
above definitions we extend to $r<0$ as $f_{\tau;r}(n)=0$ for any
$\tau$.

Our present aim is to count the number of permutations avoiding
$12\mn 3$ and avoiding (or containing exactly $r$ times) an
arbitrary generalized pattern $\tau$, since that we introduce
another notation. Let $f_{\tau;r}(n; i_1,i_2,\ldots,i_j)$ be the
number of permutations $\pi\in S_n(12\mn3)$ containing $\tau$
exactly $r$ times such that $\pi_1\pi_2\dots\pi_j=i_1i_2\dots
i_j$.

The main body of this section is divided into $6$ subsections
corresponding to the cases: $\tau$ is a general generalized
pattern, $13\mn2$, $21\mn3$, $23\mn1$, $31\mn2$, and $32\mn1$.

\subsection{$\tau$ a generalized pattern of length $k$} Here we
study certain cases of $\tau$, where $\tau$ is a generalized
pattern of length $k$ without dashes, or with exactly one dash.

\begin{theorem}\label{th11}
Let $k\geq 2$, $P_k(x)=\sum\limits_{j=0}^{k-2}\frac{x^j}{j!}$,
$G_0(x)=e^x-P_k(x)$, and for $s\geq 1$ let
$G_s(x)=G_0(x)\int_0^xG_{s-1}(t)dt$. Then
$$\mathcal{F}_{(k-1)\dots21k;r}(x)=e^{\int P_k(t)dt}\int G_k(t)dt.$$
\end{theorem}
\begin{proof}
Let $\alpha\in S_n(12\mn3)$ such that $\alpha_j=n$; so
$\alpha_1>\alpha_2>\dots>\alpha_{j-1}$. Therefore, $\alpha$
contains $\tau=(k-1)\dots21k$ exactly $r$ times if and only if
$(\alpha_{j+1},\dots,\alpha_n)$ contains $\tau$ exactly $r$ times
if $j\leq k-1$, and contains $\tau$ exactly $r-1$ times if $j\geq
k$. Thus
$$f_{\tau;r}(n)=\sum\limits_{j=1}^{k-1}\binom{n-1}{j-1}f_{\tau;r}(n-j)+\sum\limits_{j=k}^n\binom{n-1}{j-1}f_{\tau;r-1}(n-j).$$
If multiplying by $x^n/(n-1)!$ and summing over all $n\geq 1$ we
get
$$\mathcal{F}_{(k-1)\dots21k;r}(x)=\sum\limits_{j=0}^{k-2} \frac{x^j}{j!}(\mathcal{F}_{(k-1)\dots21k;r}(x)-\mathcal{F}_{(k-1)\dots21k;r-1})+e^x\mathcal{F}_{(k-1)\dots21k;r-1}(x).$$
The rest is easy to check.
\end{proof}

\begin{example}{\rm (see Claesson \cite{C})}\label{ex11}
Theorem~\ref{th11} yields for $r=0$ that
$$\mathcal{F}_{(k-1)\dots21k;0}(x)=e^{\frac{x^1}{1!}+\frac{x^2}{2!}+\dots+\frac{x^{k-1}}{(k-1)!}}.$$
If $k\rightarrow\infty$, then we get the exponential generating
function for the number of $12\mn3$-avoiding permutations in $S_n$
is given by $e^{e^x-1}$. Besides, Theorem~\ref{th11} yields for
given $k\geq 2$ and $r\rightarrow\infty$ that $G_s(x)\rightarrow
e^{e^x-\int P_k(x)}$, so
$\mathcal{F}_{(k-1)\dots21k;r}(x)\rightarrow e^{e^x-1}$.
\end{example}

Claesson~\cite{C} (see also~\cite[Pro. 28]{CM1}) proved the number
of permutations in $S_n(12\mn3,21\mn3)$ is the same number of
involutions in $S_n$. The case of varying $k$ is more interesting.
As an extension of this result.

\begin{theorem}\label{th12}
For $k\geq 2$,
$$\mathcal{F}_{(k-1)\dots21\mn k;0}(x)=e^{\frac{x^1}{1!}+\dots+\frac{x^{k-1}}{(k-1)!}}.$$
\end{theorem}
\begin{proof}
Let $\tau=(k-1)\dots 21\mn k$, by definitions we get
$$f_{\tau;0}(n)=\sum\limits_{j=1}^n f_{\tau;0}(n;j),\quad f_{\tau;0}(n;n)=f_{\tau;0}(n-1),\eqno(1)$$
and
$$f_{\tau;0}(n;i_1,\dots,i_j)=\sum\limits_{i_{j+1}=1}^{i_j-1}f_{\tau;0}(n;i_1,\dots,i_j,i_{j+1})+f_{\tau;0}(n;i_1,\dots,i_j,n)\eqno(2)$$
for all $n-1\geq i_1>i_2>\dots>i_j\geq 1$. Therefore, since
$f_{\tau;0}(n;i_1,\dots,i_j)=0$ for all $n-1>i_1>\dots>i_j\geq1$
where $j\geq k-1$ and since
$f_{\tau;0}(n;i_1,\dots,i_j,n)=f_{\tau;0}(n-1-j)$ for all
$n-1>i_1>\dots>i_j\geq 1$ where $0\leq j\leq k-2$ we have for all
$n\geq 1$
$$f_{\tau;0}(n)=\sum\limits_{j=0}^{k-2}\binom{n-1}{j}f_{\tau;0}(n-1-j).$$
The rest is easy to see as proof of Theorem \ref{th11}.
\end{proof}

In view of Example \ref{ex11} and Theorem \ref{th12} we get the
number of permutations in $S_n(12\mn3, (k-1)\dots21\mn k)$ is the
same number of permutations in $S_n(12\mn3, (k-1)\dots21k)$. In
addition,
$$S_n(12\mn3,(k-1)\dots21\mn k)=S_n(12\mn 3,(k-1)\dots21k),\eqno(3)$$
which can prove as follows. Let $\alpha=(\alpha',n,\alpha'')$;
since $\alpha$ avoids $12\mn3$ we get $\alpha'$ decreasing, so by
the principle of induction on length of $\alpha$ it is easy to see
that $\alpha$ avoids $(k-1)\dots21k$ if and only if avoids
$(k-1)\dots21\mn k$.\\

In~\cite{CM1} showed the number of permutations in $S_n(12\mn3,
13\mn2)$ is given by $2^{n-1}$. The case of varying $k$ is more
interesting. As an extension of this result.

\begin{theorem}\label{th13}
Let $k\geq 3$, then for all $n\geq 1$
$$\begin{array}{ll}
f_{(k-2)\dots21k\mn(k-1);0}(n)&=\sum\limits_{j=0}^{k-3}\binom{n-1}{j}f_{(k-2)\dots21k\mn(k-1);0}(n-1-j)+\\
&\quad\quad+\sum\limits_{j=k-2}^{n-1}\binom{n-j+k-4}{k-3}f_{(k-2)\dots21k\mn(k-1);0}(n-1-j).
\end{array}$$
\end{theorem}
\begin{proof}
Let $\tau=(k-2)\dots21k\mn(k-1)$ and let $n-1\geq
i_1>\dots>i_j\geq 1$, so for $j\leq k-3$
$$f_{\tau;0}(n;i_1,\dots,i_j,n)=f_{\tau;0}(n-1-j)$$
and for $j\geq k-2$
$$\begin{array}{ll}
f_{\tau;0}(n;i_1,\dots,i_j,n)&=f_{\tau;0}(n;n-1,n-2,\dots,n-(j-k+3),i_{j-k+4},\dots,i_j,n)=\\
&=f_{\tau;0}(n-j). \end{array}$$ Therefore, Equation (1) and
Equation (2) yield the desired result.
\end{proof}

\begin{example}{\rm(see Claesson and Mansour~\cite{CM2})}
Theorem~\ref{th13} yields for $k=3$ the number of permutations in
$S_n(12\mn3,13\mn2)$ is given by $2^{n-1}$. Another example, for
$k=4$ we get
    $$\mathcal{F}_{214\mn3;0}(x)=1+\int_0^xe^{2x+x^2/2}dx.$$
\end{example}

As a remark, similarly as proof of Equation~(3), we have for all
$k\geq3$
$$\mathcal{F}_{(k-2)\dots21k\mn(k-1);0}(x)=\mathcal{F}_{(k-2)\dots21\mn k\mn(k-1);0}(x).$$

\subsection{$\tau=13\mn2$}
\begin{theorem}\label{th14}
  Let $r$ be a nonnegative integer. Then
$$F_{13\mn2;r}(x)=\frac{1-x}{1-2x}\delta_{r,0}+\frac{x^2}{1-2x}\sum\limits_{d=1}^r
\frac{F_{13\mn2;r-d}(x)-\sum\limits_{j=0}^{d-1}f_{13\mn2;r-d}(j)x^j}{(1-x)^d}.$$
\end{theorem}
\begin{proof}
  Let $r\geq 0$, $b_r(n)=f_{13\mn2;r}(n)$, and let $1\leq i\neq j\leq n$. If $i<j$, then since
  the permutations avoiding $12\mn3$ we have $b_r(n;i,j)=0$ for
  $j\leq n-1$ and $b_r(n;i,n)=b_{r-(n-1-i)}(n-2)$, hence
  $$\sum\limits_{j=i+1}^n b_r(n;i,j) = b_{r-(n-1-i)}(n-1;n-1)=b_{r-(n-1-i)}(n-2).$$
  If $i>j$ then by definitions we have
                $$b_r(n;i,j) = b_r(n-1;j).$$
  Owing to Equation~(2) we have showed that, for
  all $1\leq i \leq n-2$,
  $$b_r(n;i)=b_{r-(n-1-i)}(n-2)+\sum\limits_{j=1}^{i-1}b_r(n-1;j).\eqno(4)
  $$
  Moreover, it is plain that
  $$b_r(n;n)=b_r(n;n-1)=b_r(n-1),\eqno(5)$$
  and by means of induction we shall show that Equation~(4)
  implies: If $2\leq m \leq n-1$ then
  $$\begin{array}{l}
    b_r(n;n-m)=\sum\limits_j (-1)^j\left[
          \binom{m-1}{j}+\binom{m-2}{j-1}
        \right]b_r(n-1-j)+\\
         \qquad\qquad\qquad+b_{r-(m-1)}(n-2)-
        \sum\limits_{d\geq 1}\sum\limits_j(-1)^j\binom{m-2-d}{j}b_{r-d}(n-3-j).
      \end{array}\eqno(6)$$
  First we verify the statement for
  $m=2$; in this case Equation~(6) becomes
  $$b_r(n;n-2) = b_r(n-1) - 2 b_r(n-2) + b_{r-1}(n-2).$$
  Indeed,
  $$\begin{array}{l}
     b_r(n;n-m) =\\
     \quad\quad = \sum\limits_{j=1}^{n-m-1} b_r(n;n-m,j)+ \sum\limits_{j=n-m+1}^{m}b_r(n;n-m,j)\\
     \quad\quad = \sum\limits_{j=1}^{n-m-1} b_r(n-1;j)+ b_{r-(m-1)}(n-1;n-1)\\
      \quad\quad = b_r(n-1) - 2 b_r(n-2) + b_{r-(m-1)}(n-2) -
      \sum\limits_{k=2}^{m-1} b_r(n-1;n-1-k),
    \end{array}\eqno(7)$$
  where the three equalities follow from Equation~(2),
  and Equation~(1) together with ~(4) and ~(5),
  respectively. Now simply put $m=2$ to obtain Equation~(6).

  Assume that Equation~(7) holds for all $k$ such that $2\leq k
  \leq m-1$. Then, employing the familiar identity $\binom 1 k +
  \binom 2 k + \cdots +\binom n k = \binom{n+1}{k+1}$, the trailing
  sum in Equation~(7) expands as follows. Since
$$\begin{array}{l}
\sum\limits_{k=2}^{m-1}\sum\limits_j(-1)^j\Biggl[\binom{k-1}{j}+\binom{k-2}{j-1}\Biggr]
b_r(n-2-j)=\\
\quad=\sum\limits_{j}(-1)^j\Biggl[\binom{m-1}{j+1}+\binom{k-2}{j}\Biggr]
b_r(n-2-j) -b_r(n-2)=\\
\quad\quad=-\sum\limits_{j}(-1)^j\Biggl[\binom{k-1}{j}+\binom{k-2}{j-1}\Biggr]b_r(n-1-j)
+b_r(n-1)-2b_r(n-2),
\end{array}$$

$$\sum\limits_{k=2}^{m-1}b_{r-(k-1)}(n-3)=\sum\limits_{d\geq1}b_{r-d}(n-3),$$
and
$$\begin{array}{l}
\sum\limits_{k=2}^{m-1}\sum\limits_{d\geq1}\sum\limits_j(-1)^j\binom{k-2-d}{j}b_{r-d}(n-4-j)=\\
\qquad\qquad\quad=-\sum\limits_{d\geq
1}\sum\limits_j(-1)^j\binom{m-2-d}{j}b_{r-d}(n-3-j)
+\sum\limits_{d\geq 1}b_{r-d}(n-3)
\end{array}$$
with using Equation~(7) we get that Equation~(6) holds for $k=m$,
by the principle of induction the universal validity of
Equation~(6) follows.\\

Now, if summing  $b_r(n;n-m)$ over all $0\leq m\leq n-1$, then by
using Equation~(1), ~(5), and ~(6) we get
$$\begin{array}{l}
\sum\limits_j(-1)^j\left[\binom{n-1}{j}+\binom{n-2}{j-1}\right]b_r(n-j)=\\
\qquad\qquad\qquad\qquad=\sum\limits_{d\geq1}\sum\limits_j(-1)^j\binom{n-2-d}{j}b_{r-d}(n-2-j).
\end{array}\eqno(8)$$

Using \cite[Lem.~7]{CM2} to transfer the above equation in terms
of ordinary generating functions
$$\begin{array}{l}
(1-u)F_{13\mn2;r}\left(\frac{u}{1+u}\right)=\delta_{r,0}+\\
\quad\quad\quad+u^2\sum\limits_{d=1}^r (1+u)^{d-1}\left[
F_{13\mn2;r-d}\left(\frac{u}{1+u}\right)-\sum\limits_{j=0}^{d-1}f_{13\mn2;r-d}(j)\left(\frac{u}{1+u}\right)^j\right].
\end{array}$$
Putting $x=u/(1+u)$ ($u=x/(1-x)$) we get the desired result.
\end{proof}

An application for Theorem \ref{th14} we get the exact formula for
$f_{13\mn2;r}(n)$ for $r=0,1,2,3,4$.

\begin{corollary}
For all $n\geq 1$,
$$\begin{array}{l}
f_{13\mn2;0}(n)=2^{n-1};\\
f_{13\mn2;1}(n)=(n-3)2^{n-2}+1;\\
f_{13\mn2;2}(n)=(n^2-3n -6)2^{n-4}+n;\\
f_{13\mn2;3}(n)=1/3(n^3-31n-18)2^{n-5}+n^2-n+1;\\
f_{13\mn2;4}(n)=1/3(n-1)(n^3+7n^2-546n-312)2^{n-8}+2/3(n-1)(n^2-2n+3).
\end{array}$$
\end{corollary}
\subsection{$\tau=21\mn3$}
By definitions it is easy to obtain the following:

\begin{lemma}\label{l161}
Let $n\geq 1$; then
$$\begin{array}{l}
f_{21\mn3;r}(n)=f_{21\mn3;r}+\sum\limits_{i=1}^{n-1}f_{21\mn3;r}(n;i),\\
f_{21\mn3;r}(n;i)=f_{21\mn3;r}(n-2)+\sum\limits_{j=1}^{i-1}f_{21\mn3;r-(n-i)}(n-2;j),\quad\mbox{for}
\ 1\leq i\leq n-1. \end{array}$$
\end{lemma}

Using the above lemma for given $r$ we obtain the exact formula
for $f_{21\mn3;r}$. Here we present the first three cases
$r=0,1,2$.

\begin{theorem}\label{th16}
For all $n\geq 1$
$$\begin{array}{l}
f_{21\mn3;0}(n)=f_{21\mn3;0}(n-1)+(n-1)f_{21\mn3;0}(n-2);\\
\ \\
f_{21\mn3;1}(n)=f_{21\mn3;1}(n-1)+(n-1)f_{21\mn3;1}(n-2)+f_{21\mn3;0}(n-1)-f_{21\mn3;0}(n-2);\\
\ \\
f_{21\mn3;2}(n)=f_{21\mn3;2}(n-1)+(n-1)f_{21\mn3;2}(n-2)+f_{21\mn3;1}(n-1)-f_{21\mn3;1}(n-2)+\\
\qquad\qquad\qquad\qquad\quad\qquad\qquad\qquad\qquad\qquad\qquad+f_{21\mn3;0}(n-1)-2f_{21\mn3;0}(n-2).
\end{array}$$
\end{theorem}
\begin{proof}
Case $r=0$: Lemma~\ref{l161} yields
$f_{21\mn3;0}(n;i)=f_{21\mn3;0}(n-1)$ for all $1\leq i\leq n-1$,
so for all $n\geq 1$
$$f_{21\mn3;0}(n)=f_{21\mn3;0}(n-1)+(n-1)f_{21\mn3;0}(n-2).$$

Case $r=1$: Lemma~\ref{l161} yields
$f_{21\mn3;1}(n;i)=f_{21\mn3;1}(n-2)$ for all $1\leq i\leq n-2$,
and
$$f_{21\mn3;1}(n;n-1)=\sum\limits_{j=1}^{n-2}f_{21\mn3;0}(n-1;j)+f_{21\mn3;1}(n-2)$$
which equivalent to (by use the case $r=0$)
$$f_{21\mn3;1}(n;n-1)=f_{21\mn3;1}(n-2)+f_{21\mn3;0}(n-1)-f_{21\mn3;0}(n-2).$$
Therefore, for all $n\geq 1$
$$f_{21\mn3;1}(n)=f_{21\mn3;1}(n-1)+(n-1)f_{21\mn3;1}(n-2)+f_{21\mn3;0}(n-1)-f_{21\mn3;0}(n-2).$$

Case $r=2$: Similarly as the cases $r=0,1$.
\end{proof}
\subsection{$\tau=23\mn1$}

\begin{theorem}\label{th17}
For any $r$ nonnegative integer,
$$\begin{array}{l}
F_{23\mn1;r}(x)=\frac{\delta_{r,0}}{1-x}+\\
\qquad\qquad\qquad+x^2\sum\limits_{d=0}^r(1-x)^{j-2}\left[F_{23\mn1;r-d}
\left(\frac{x}{1-x}\right)-\sum\limits_{j=0}^{d-1}f_{23\mn1;r-d}(j)\left(\frac{x}{1-x}\right)^j\right].
\end{array}$$
\end{theorem}
\begin{proof}
By definition it is easy to state the following statement:

\begin{lemma}\label{l171}
Let $n\geq 1$; then
$$\begin{array}{l}
f_{23\mn1;r}(n)=f_{23\mn1;r}(n-1)+\sum\limits_{i=1}^{n-1}f_{23\mn1;r}(n;i),\\
f_{23\mn1;r}(n;i)=f_{23\mn1;r-(i-1)}(n-2)+\sum\limits_{j=1}^{i-1}f_{23\mn1;r}(n-1;j),\quad
\mbox{for}\ 1\leq i\leq n-1. \end{array}$$
\end{lemma}

By consider the same argument proof of Equation~(6) with use
Lemma~\ref{l171} we have for all $1\leq m\leq n-1$,
$$f_{23\mn1;r}(n;m)=f_{23\mn1;r+1-m}(n-2)+\sum\limits_{d\geq1}\sum\limits_j\binom{m-1-d}{j}f_{23\mn1;r+1-d}(n-3-j).$$
Therefore, by summing $f_{23\mn1;r}(n;m)$ over all $1\leq m\leq n$
we shall show that Lemma~\ref{l171} implies for all $n\geq 1$,
$$f_{23\mn1;r}(n)=f_{23\mn1;r}(n-1)+\sum\limits_{d=0}^r\sum\limits_{j=0}^{n-2-d}\binom{n-2-d}{j}f_{23\mn1;r-d}(n-2-j).$$
Hence, using \cite[Lem.~7]{CM2} we get the desired result.
\end{proof}

\begin{example}{\rm (see Claesson and Mansour
\cite[Pro.~7]{CM1})}\label{exxx1} Theorem~\ref{th17} for $r=0$
yields
$$F_{23\mn1;0}(x)=\frac{1}{1-x}+\frac{x^2}{(1-x)^2}F_{23\mn1;0}\left(\frac{x}{1-x}\right).$$
An infinite number of application of this identity we have
$$F_{23\mn1;0}(x)=\sum\limits_{k\geq0}\frac{x^{2k}}{p_{k-1}(x)p_{k+1}(x)},$$
where $p_m(x)=\prod_{j=0}^m(1-jx)$. An another example,
Theorem~\ref{th17} for $r=1$ yields (similarly)
$$F_{23\mn1;1}(x)=\sum\limits_{d\geq 0}\left[ \frac{x^{2d+2}}{p_{d+1}(x)}\left(
\sum\limits_{k\geq0}\frac{x^{2k}}{p_{k+d-1}(x)p_{k+d+1}(x)}-1\right)\right].$$
\end{example}
\subsection{$\tau=31\mn2$}
By definitions it is easy to state the following:

\begin{lemma}\label{l151}
Let $n\geq 1$; then
$$\begin{array}{l}
f_{31\mn2;r}(n;n)=\sum\limits_{j=1}^{n-1}f_{31\mn2;r+1-j}(n-1;n-j),\\
f_{31\mn2;r}(n;i)=f_{31\mn2;r}(n-1;n-1)+\sum\limits_{j=1}^{i-1}f_{31\mn2;r-(i-1-j)}(n-1;j),\quad\mbox{for}\
1\leq i\leq n-1.
\end{array}$$
\end{lemma}

\begin{theorem}\label{th15}
  Let $r$ be a nonnegative integer, then $f_{31\mn2;r}(n)$ is a
  polynomial of degree at most $2r+2$ with coefficient in $\mathcal{Q}$, where $n\geq
  0$.
\end{theorem}
\begin{proof}
Using Lemma~\ref{l151} for $r=0$ we obtain that, first
$f_{31\mn2;0}(n;n)=1$ and $f_{31\mn2;0}(n;1)=1$, second
$f_{31\mn2;0}(n;j)=j$ for all $1\leq j\leq n-1$. Hence, for all
$n\geq 0$,
$$f_{31\mn2;0}(n)=\binom{n}{2}+1.$$

Now, assume that $f_{31\mn2;d}(n;j)$ is a polynomial of degree at
most $2d+1$ with coefficient in $\mathcal{Q}$ for all $1\leq j\leq
n$ where $d=0,1,2,\dots,r-1$. Therefore, Lemma~\ref{l151} with
induction hypothesis imply, first $f_{31\mn2;r}(n;n)$ and
$f_{31\mn2;r}(n;1)$ are polynomials of degree at most $2r$, and
then $f_{31\mn2;r}(n;j)$ is a polynomial of degree at most $2r+1$.
So, by use the principle of induction on $r$ we get that
$f_{31\mn2;r}(n;j)$ is a polynomial of degree at most $2r+1$ with
coefficient in $\mathcal{Q}$ for all $r\geq 0$. Hence, since
$f_{31\mn2;r}(n)=\sum\limits_{j=1}^n f_{31\mn2;r}(n;j)$ we get the
desired result.
\end{proof}

An application for Theorem ~\ref{th15} with the initial values of
the sequence $f_{31\mn2;r}(n)$ we have the exact formula for
$f_{31\mn2;r}(n)$ where $r=0,1,2,3$.

\begin{corollary}
For all $n\geq 0$,
$$\begin{array}{l}
f_{31\mn2;0}(n)=1+n(n-1)/2;\\
f_{31\mn2;1}(n)=n(n-1)(n-2)(3n-5)/24;\\
f_{31\mn2;2}(n)=n(n-1)(n-2)(n-3)(5n^2-3n-38)/720;\\
f_{31\mn2;3}(n)=n(n-1)(n-2)(n-3)(n-4)(7n^3+10n^2+205n-1142)/40320.
\end{array}$$
\end{corollary}

\subsection{$\tau=32\mn1$}
By definitions it is easy to state the following:

\begin{lemma}\label{l152}
Let $n\geq 1$; then
$$\begin{array}{l}
f_{32\mn1;r}(n;n)=\sum\limits_{j=1}^{n-1}f_{31\mn2;r+1-j}(n-1;j),\\
f_{32\mn1;r}(n;i)=f_{32\mn1;r}(n-1;n-1)+\sum\limits_{j=1}^{i-1}f_{31\mn2;r+1-j)}(n-1;j),\quad\mbox{for}\
1\leq i\leq n-1.
\end{array}$$
\end{lemma}

\begin{theorem}\label{th15b}
  Let $r$ be a nonnegative integer, then $f_{32\mn1;r}(n)$ is a
  polynomial of degree at most $r+1$ with coefficient in $\mathcal{Q}$,
  where $n\geq r+2$.
\end{theorem}
\begin{proof}
Using Lemma~\ref{l152} for $r=0$ we obtain that, first
$f_{32\mn1;0}(n;n)=1$ and $f_{32\mn1;0}(n;1)=1$, second
$f_{31\mn2;0}(n;j)=2$ for all $2\leq j\leq n-1$. Hence, for all
$n\geq 2$,
                $$f_{31\mn2;0}(n)=2n-2.$$

Let $n\geq r+2$ and let us assume that $f_{32\mn1;d}(n;j)$ is a
polynomial of degree at most $d$ with coefficient in $\mathcal{Q}$
for all $1\leq j\leq n$ where $d=0,1,2,\dots,r-1$.
Lemma~\ref{l152} with induction hypothesis imply, first
$f_{32\mn1;r}(n;n)$ and $f_{32\mn1;r}(n;1)$ are polynomials of
degree at most $r$ with coefficient in $\mathcal{Q}$, and then
$f_{31\mn2;r}(n;j)$ where $2\leq j\leq n-1$ is a polynomial of
degree at most $r$ with coefficient in $\mathcal{Q}$. So, by use
the principle of induction on $r$ we get that $f_{32\mn1;r}(n;j)$
is a polynomial of degree at most $r$ with coefficient in
$\mathcal{Q}$ for all $r\geq 0$. Hence, since
$f_{32\mn1;r}(n)=\sum\limits_{j=1}^n f_{32\mn1;r}(n;j)$ we get the
desired result.
\end{proof}

An application for Theorem ~\ref{th15b} with the initial values of
the sequence $f_{32\mn1;r}(n)$ we get the exact formula for
$f_{32\mn1;r}(n)$ for $r=0,1,2,3$.

\begin{corollary}
$$\begin{array}{ll}
\mbox{For all}\ n\geq 2,\ & f_{32\mn1;0}(n)=2n-2;\\
\mbox{For all}\ n\geq 3,\ & f_{32\mn1;1}(n)=(n-3)(2n-1);\\
\mbox{For all}\ n\geq 4,\ & f_{32\mn1;2}(n)=(n-4)(n^2-3n+1);\\
\mbox{For all}\ n\geq 5,\ &
f_{32\mn1;3}(n)=(n-5)(2n^3-13n^2+47n-6)/6.
\end{array}$$
\end{corollary}
\section{Avoiding $13\mn2$}
Let $g_{\tau;r}(n)$ be the number of all permutations in
$S_n(13\mn2)$ containing $\tau$ exactly $r$ times. The
corresponding ordinary generating function we denote by
$G_{\tau;r}(x)$; that is,
$G_{\tau;r}(x)=\sum\limits_{n\geq0}g_{\tau;r}(n)x^n$. The above
definitions we extend to $r<0$ as $g_{\tau;r}(n)=0$ for any
$\tau$.

In the current section, our present aim is to count the number of
permutations avoiding $13\mn2$ and containing $\tau$ exactly $r$
times where $\tau$ a generalized pattern of type $(2,1)$, and
since that we introduce another notation. Let $g_{\tau;r}(n;
i_1,i_2,\ldots,i_j)$ be the number of permutations $\pi\in
S_n(13\mn2)$ containing $\tau$ exactly $r$ times such that
$\pi_1\pi_2\dots\pi_j=i_1i_2\dots i_j$.

The main body of this section is divided to three subsections
corresponding to the cases $\tau$ is $12\mn3$; $21\mn3$; and
$23\mn1$, or $31\mn2$, or $32\mn1$.

\subsection{$\tau=12\mn3$}
\begin{theorem}\label{th21}
Let $r$ be any nonnegative integer number; then there exists a
polynomials $p_r(n)$ and $q_{r-1}(n)$ of degree at most $r$ and
$r-1$ respectively, with coefficient in $\mathcal{Q}$ such that
for all $n\geq 1$,
$$g_{12\mn3;r}(n)=p_r(n)\cdot 2^n+q_{r-1}(n).$$
\end{theorem}
\begin{proof}
Let $r\geq 0$, and let $1\leq i\neq j\leq n$. If $i<j$, then since
the permutations avoiding $13\mn2$ we have $g_{12\mn3;r}(n;i,j)=0$
for $i+2\leq j\leq n$ and
$g_{12\mn3;r}(n;i,i+1)=g_{12\mn3;r-(n-1-i)}(n-1;i)$, hence
        $$\sum\limits_{j=i+1}^n g_{12\mn3;r}(n;i,j)=g_{12\mn3;r-(n-1-i)}(n-1;i).$$
If $i>j$ then by definitions we have
                $$g_{12\mn3;r}(n;i,j)= g_{12\mn3;r}(n-1;j).$$
Owing to the definitions we have showed that, for all $1\leq i\leq
n-2$,
$$g_{12\mn3;r}(n;i)=g_{12\mn3;r-(n-1-i)}(n-1;i)+\sum\limits_{j=1}^{i-1}g_{12\mn3;r}(n-1;j).\eqno(1')$$
Moreover, it is plain that
  $$g_{12\mn3;r}(n;n)=g_{12\mn3;r}(n;n-1)=g_{12\mn3;r}(n-1),\eqno(2')$$
and for all $1\leq j\leq n-r-2$
    $$g_{12\mn3;r}(n;j)=0.\eqno(3')$$

Now we ready to prove the theorem. Let $r=0$; by Equation~(3') we
get $g_{12\mn3;0}(n;j)=0$ for all $j\leq n-2$ and by Equation~(2')
we have  $g_{12\mn3;0}(n;n-1)=g_{12\mn3;0}(n;n)=d_0(n-1)$, so
$d_r(n)=2^{n-1}$. Therefore, the theorem holds for $r=0$.

Let $r\geq 1$, and let us assume that for all $0\leq m\leq s-1$
and all $0\leq s\leq r-1$ there exists a polynomials $p_{m}(n)$
and $q_{m-1}(n)$ of degree at most $m$ and $m-1$ respectively with
coefficient in $\mathcal{Q}$ such that
$g_{12\mn3;s}(n;n-s-1+m)=p_{m}(n)2^n+q_{m-1}(n)$, and there exists
a polynomials $v_s(n)$ and $u_{s-1}(n)$ of degree at most $s$ and
$s-1$ respectively with coefficient in $\mathcal{Q}$ such that
$g_{12\mn3;s}(n-m)=v_s(n)2^n+u_{s-1}(n)$ where $m=0,1$.

So, using Equation~(1') for $m=0,1,\dots,r-1$ and the induction
hypothesis imply that there exists a polynomials $a_{m}(n)$ and
$b_{m-1}(n)$ of degree at most $m$ and $m-1$ respectively with
coefficient in $\mathcal{Q}$ such that
$$g_{12\mn3;r}(n;n-r-1+m)=a_{m}(n)2^n+b_{m-1}(n).$$
Besides, Owing to Equations~(1'), (2'), and ~(3') we have showed
that
    $$g_{12\mn3;r}(n)=2g_{12\mn3;r}(n-1)+\sum\limits_{j=2}^{r+1}g_{12\mn3;r}(n;n-j),$$
which means that $g_{12\mn3;r}(n)$ is given by
$v_r(n)2^n+u_{r-1}(n)$ and
$g_{12\mn3;r}(n;n)=g_{12\mn3;r}(n;n-1)=g_{12\mn3;r}(n-1)$.
Therefore, the statement holds for $s=r$. Hence, by the principle
of induction on $r$ the theorem holds.
\end{proof}

An application for Theorem \ref{th21} with the initial values of
the sequence $g_{12\mn3;r}(n)$ we obtain the exact formula for
$g_{12\mn3;r}(n)$ for $r=0,1,2,3$.

\begin{corollary}
For all $n\geq 1$;
$$\begin{array}{l}
g_{12\mn3;0}(n)=2^{n-1};\\
g_{12\mn3;2}(n)=(n-3)2^{n-2}+1;\\
g_{12\mn3;2}(n)=(n^2-11n+34)2^{n-4}-n-2;\\
g_{12\mn3;3}(n)=1/3(n^3-24n^2+257n-954)2^{n-5}+n^2+4n+10.
\end{array}$$
\end{corollary}
\subsection{$\tau=21\mn3$}
\begin{theorem}\label{th22}
Let $r$ be any nonnegative integer number. Then there exist a
polynomial $p_r(n)$ of degree at most $r$ with coefficient in
$\mathcal{Q}$, such that for all $n\geq r$
  $$g_{21\mn3;r}(n)=p_r(n)\cdot 2^n.$$
\end{theorem}
\begin{proof}
By definitions it is easy to state

\begin{lemma}\label{l221}
Let $n\geq 1$; then
$$\begin{array}{l}
g_{21\mn3;r}(n)=g_{21\mn3;r}(n-1)+\sum\limits_{i=1}^{n-1}g_{23\mn1;r}(n;i),\\
g_{21\mn3;r}(n;i)=g_{21\mn3;r}(n-1;i)+\sum\limits_{j=1}^{i-1}g_{21\mn3;r-(n-i)}(n-1;j),\quad\mbox{for}\
1\leq i\leq n-1.
\end{array}$$
\end{lemma}

Lemma~\ref{l221} implies for $r=0$ as follows. First
$g_{21\mn3;0}(n;m)=2^{m-2}$ for all $m\geq 2$, and second
$g_{21\mn3;0}(n;1)=1$. Hence $g_{21\mn3;0}(n)=2^{n-1}$, so the
theorem holds for $r=0$.

Let $r\geq1$ and assume that for $2\leq m\leq n-1$ the expression
$\sum\limits_{j=1}^m g_{21\mn3;d}(n;j)$ is given by $q_d^m(n)2^m$
where $q_{d}^m(n)$ is a polynomial of degree at most $d$ with
coefficient in $\mathcal{Q}$ for all $0\leq d\leq r-1$.
Lemma~\ref{l221} yields
$$\begin{array}{ll}
g_{21\mn3;r}(n;1)&=g_{21\mn3;r}(n-1;1),\\
g_{21\mn3;r}(n;2)&=g_{21\mn3;r}(n-1;2),\\
\qquad\qquad\qquad\qquad\vdots\\
g_{21\mn3;r}(n;n-r-1)&=g_{21\mn3;r}(n-1;n-r-1),\\
g_{21\mn3;r}(n;n-r+1)&=g_{21\mn3;r}(n-1;n-r+1)+\sum\limits_{j=1}^{n-r}g_{21\mn3;1}(n-1;j),\\
\qquad\qquad\qquad\qquad\vdots\\
g_{21\mn3;r}(n;n-1)&=g_{21\mn3;r}(n-1;n-1)+\sum\limits_{j=1}^{n-2}g_{21\mn3;r-1}(n-1;j),\\
g_{21\mn3;r}(n;n)&=g_{21\mn3;r}(n-1),
\end{array}$$
with induction hypothesis imply for $2\leq m\leq n-1$
$$\sum\limits_{j=1}^m g_{21\mn3;r}(n;j)=\sum\limits_{j=1}^mg_{21\mn3;r}(n-1;j)+q_{r-1}^m(n)2^m,$$
where $q_{r-1}^m(n)$ is a polynomial of degree at most $r-1$ with
coefficient in $\mathcal{Q}$. Therefore, for $2\leq m\leq n-1$
$\sum\limits_{j=1}^mg_{12\mn3;r}(n;j)$ can be expressed by
$q_r^m(n)2^m$ where $q_r^m(n)$ is a polynomial of degree at most
$r$ with coefficient in $\mathcal{Q}$. Hence, with using
Lemma~\ref{l221} we get there exist a polynomial $a_r(n)$ of
degree at most $r$ with coefficient in $\mathcal{Q}$ such that
$$g_{12\mn3;r}(n)=g_{12\mn3;r}(n-1)+a_r(n)2^n,$$
so the theorem holds.
\end{proof}

Using Theorem \ref{th22} with the initial values of the sequence
$g_{21\mn3;r}(n)$ for $r=0,1,2,3$ we get

\begin{corollary}\

$\begin{array}{lll}
{\rm (i)} &\mbox{For all}\ n\geq1, &g_{21\mn3;0}(n)=2^{n-1};\\

{\rm (ii)} &\mbox{For all}\ n\geq2, &g_{21\mn3;1}(n)=(n-2)2^{n-3};\\

{\rm (iii)} &\mbox{For all}\ n\geq3,&
g_{21\mn3;2}(n)=(n^2+n-12)2^{n-6};\\

{\rm (iv)} & \mbox{For all}\ n\geq4,&
g_{21\mn3;3}(n)=1/3(n-4)(n^2+13n+6)2^{n-8}.\\
\end{array}$
\end{corollary}

\subsection{$\tau=23\mn1$, $\tau=31\mn2$, or $\tau=32\mn1$}
Similarly, using the argument proof of Theorem~\ref{th22} with the
principle of induction yield

\begin{theorem}\label{th23}
Let $r$ be any nonnegative integer number. Then

{\rm (i)} there exists a polynomial $p_{r-1}(n)$ of degree at most
$r-1$ with coefficient in $\mathcal{Q}$ and a constant $c$, such
that for all $n\geq r$
  $$g_{23\mn1;r}(n)=c\cdot2^n+p_{r-1}(n).$$

{\rm (ii)} there exists a polynomials $p_r(n)$ and $q_{2r-2}(n)$
of degree at most $r$ and $2r-2$ respectively; with coefficient in
$\mathcal{Q}$ such that for all $n\geq 1$
    $$g_{31\mn2;r}(n)=p_r(n)2^n+q_{2r-2}(n).$$

{\rm (iii)} there exist a polynomial $p_{r+2}(n)$ of degree at
most $r+2$ with coefficient in $\mathcal{Q}$ such that for all
$n\geq r$
  $$g_{32\mn1;r}(n)=p_{r+2}(n).$$
\end{theorem}

Using Theorem \ref{th23} with the initial values of the sequences
$g_{23\mn1;r}(n)$, $g_{31\mn2;r}(n)$ and $g_{32\mn1;r}(n)$ for
$r=0,1,2,3,4$ we get the following:

\begin{corollary}\

$\begin{array}{lll}
{\rm (i)}  & \mbox{For all}\ n\geq1,& g_{23\mn1;0}(n)=2^{n-1};\\
{\rm (ii)} & \mbox{For all}\ n\geq2,& g_{23\mn1;1}(n)=2^{n-2}-1;\\
{\rm (iii)}& \mbox{For all}\ n\geq3,& g_{23\mn1;2}(n)=2^{n-1}-n-1;\\
{\rm (iv)} & \mbox{For all}\ n\geq4,& g_{23\mn1;3}(n)=5\cdot2^{n-3}-1/2(n^2-n+8);\\
{\rm (v)}  & \mbox{For all}\ n\geq5,&
g_{23\mn1;4}(n)=2^{n}-1/6(n+1)(n^2-4n+24).
\end{array}$
\end{corollary}

\begin{corollary}
For all $n\geq 1$;

$\begin{array}{ll}
{\rm (i)} & g_{31\mn2;0}(n)=2^{n-1};\\
{\rm (ii)} & g_{31\mn2;1}(n)=(n-3)2^{n-2}+1;\\
{\rm (iii)} & g_{31\mn2;2}(n)=(n^2-3n-14)2^{n-4}+1/2(n^2+n+12);\\
{\rm (iv)} &
g_{31\mn2;3}(n)=1/3(n^3-55n-90)2^{n-5}+1/12(n^4+11n^2+12n+12).
\end{array}$
\end{corollary}

\begin{corollary}
For all $n\geq 1$;

$\begin{array}{ll}
{\rm (i)}  & g_{32\mn1;0}(n)=1/2n(n-1)+1;\\
{\rm (ii)} & g_{32\mn1;1}(n)=1/6(n-1)(n-2)(2n-3);\\
{\rm (iii)} & g_{32\mn1;2}(n)=1/6(n-2)(n-3))2n-5);\\
{\rm (iv)} & g_{32\mn1;3}(n)=1/8(n-3)(n^3-3n^2-10n+32);\\
{\rm (v)} & g_{32\mn1;4}(n)=1/24(n-4)(3n^3-10n^2-55n+198).
\end{array}$
\end{corollary}

\section{Avoiding $21\mn3$}
Let $h_{\tau;r}(n)$ be the number of all permutations in
$S_n(21\mn3)$ containing $\tau$ exactly $r$ times. The
corresponding exponential and ordinary generating function are
denoted by $\mathcal{H}_{\tau;r}(x)$ and $H_{\tau;r}(x)$
respectively; that is,
$\mathcal{H}_{\tau;r}(x)=\sum\limits_{n\geq0}\frac{h_{\tau;r}(n)}{n!}x^n$
and $H_{\tau;r}(x)=\sum\limits_{n\geq0}h_{\tau;r}(n)x^n$. The
above definitions are extended to $r<0$ as $h_{\tau;r}(n)=0$ for
any $\tau$.

In the current section, our present aim is to count the number of
permutations avoiding $21\mn3$ and containing $\tau$ exactly $r$
times where $\tau$ a generalized pattern of type $(2,1)$, and
since that we introduce another notation. Let $h_{\tau;r}(n;
i_1,i_2,\ldots,i_j)$ be the number of permutations $\pi\in
S_n(21\mn3)$ containing $\tau$ exactly $r$ times such that
$\pi_1\pi_2\dots\pi_j=i_1i_2\dots i_j$.

The main body of the current section is divided to five
subsections corresponding to the cases $\tau$ is a general
generalized pattern; $12\mn3$; $13\mn2$, $31\mn2$; $23\mn1$; and
$32\mn1$.

\subsection{$\tau$ is a general generalized pattern} Here we
study certain cases of $\tau$, where $\tau$ is a generalized
pattern of length $k$ without dashes, or with exactly one dash.\\

First of all let us define a bijection $\Phi$ between the set
$S_n(12\mn3)$ and the set $S_n(21\mn3)$ as follows. Let
$\pi=(\pi',n,\pi'')$, where $n$ the maximal element of $\pi$, be
any $12\mn3$-avoiding permutation of $s$ elements; we define by
induction
            $$\Phi(\pi)=(R(\pi'),n,\Phi(\pi'')),$$
where $R(\pi')$ is the reversal of $\pi'$. Since $\pi$ is
$12\mn3$-avoiding permutation we have $\pi_1>\dots>\pi_{j-1}$ so
by using the principle of induction on length $\pi$ we get
$\Phi(\pi)$ is $21\mn3$-avoiding permutation. Also, it is easy to
see by using the principle of induction that $\Phi^{-1}=\Phi$,
hence $\Phi$ is a bijection.

\begin{theorem}\label{th31}
For all $k\geq 1$;
$$
\mathcal{H}_{12\dots(k-1)k;0}(x)=\mathcal{F}_{(k-1)\dots21k;0}(x),\qquad
\mathcal{H}_{12\dots(k-1)k;1}(x)=\mathcal{F}_{(k-1)\dots21k;1}(x).$$
\end{theorem}
\begin{proof}
Using the bijection $\Phi:S_n(12\mn3)\rightarrow S_n(21\mn3)$ we
get the desired result: the permutation $\pi\in S_n(12\mn3)$
contains $(k-1)\dots1k$ exactly $r$ ($r=0,1$) times if and only if
the permutations $\Phi(\pi)$ contains $12\dots(k-1)k$ exactly $r$.
\end{proof}

\begin{example}{\rm (see Claesson \cite{C})}\label{ex31}
Theorem~\ref{th11} and Theorem~\ref{th31} yield for $r=0$ that
$$\mathcal{H}_{12\dots k;0}(x)=e^{\frac{x^1}{1!}+\frac{x^2}{2!}+\dots+\frac{x^{k-1}}{(k-1)!}}.$$
If $k\rightarrow\infty$, then we get the exponential generating
function for the number of $21\mn3$-avoiding permutations in $S_n$
is given by $e^{e^x-1}$.
\end{example}

In \cite{C,CM1} proved the number of permutations in $S_n(21\mn3)$
avoiding $12\mn3$ is the same number of involutions in $S_n$. The
case of varying $k$ is more interesting. As an extension of these
results (the proofs are immediately holds by using the bijection
$\Phi$).

\begin{theorem}\label{th32}
For $k\geq 1$,
$$\begin{array}{l}
\mathcal{H}_{12\dots(k-1)\mn
k;0}(x)=e^{\frac{x^1}{1!}+\dots+\frac{x^{k-1}}{(k-1)!}};\\
\mathcal{H}_{12\dots(k-2)k\mn(k-1);0}(x)=\mathcal{H}_{12\dots(k-2)\mn
k\mn(k-1);0}(x) =\mathcal{F}_{(k-2)\dots21k\mn(k-1);0}(x).
\end{array}$$
\end{theorem}

Using the bijection $\Phi$ we get easily other results as follows.

\begin{theorem}
(i) The number of permutations in $S_n$ containing $12\mn3$
exactly once is the same number of permutations containing
$21\mn3$ exactly once;

(ii) The number of permutations in $S_n$ containing $12\mn3$ once
and containing $(k-1)\dots21\mn k$ (resp. $(k-1)\dots21k$) exactly
$r=0,1$ times, is the same number of permutations in $S_n$
containing $21\mn3$ once and containing $12\dots(k-1)\mn k$ (resp.
$(12\dots (k-1)k$) exactly $r=0,1$ times.
\end{theorem}
\subsection{$\tau=12\mn3$}
\begin{theorem}\label{th35}
For all $n\geq 1$,
$$\begin{array}{l}
h_{12\mn3;0}(n)=h_{12\mn3;0}(n-1)+(n-1)h_{12\mn3;0}(n-2);\\
\ \\
h_{12\mn3;1}(n)=h_{12\mn3;1}(n-1)+(n-1)h_{12\mn3;1}(n-2)+(n-2)h_{12\mn3;0}(n-3);\\
\ \\
h_{12\mn3;2}(n)=h_{12\mn3;2}(n-1)+(n-1)h_{12\mn3;2}(n-2)+(n-2)h_{12\mn3;1}(n-3)+\\
\qquad\qquad\qquad\qquad\qquad\qquad\qquad\qquad\qquad\qquad\quad\quad+(n-3)h_{12\mn3;0}(n-3).
\end{array}$$
\end{theorem}
\begin{proof}
Immediately, definitions yield the following statement:

\begin{lemma}\label{l351}
Let $n\geq 1$; then
$$\begin{array}{l}
h_{12\mn3;r}(n)=h_{12\mn3;r}(n-1)+\sum\limits_{j=1}^{n-1}h_{12\mn3;r}(n;j),\\
h_{12\mn3;r}(n;j)=h_{12\mn3;r}(n-2)+\sum\limits_{i=j}^{n-1}h_{12\mn3;r-(n-i-1)}(n-1;i),\quad\mbox{for}\
1\leq j\leq n-1.
\end{array}$$
\end{lemma}

Case $r=0$: Lemma~\ref{l351} yields
$h_{12\mn3;0}(n;j)=h_{12\mn3;0}(n-2)$ for all $1\leq j\leq n-1$,
hence
$$h_{12\mn3;0}(n)=h_{12\mn3;0}(n-1)+(n-1)h_{12\mn3;0}(n-2).$$

Case $r=1$: Lemma~\ref{l351} yields
$$h_{12\mn3;1}(n;j)=h_{12\mn3;1}(n-2)+h_{12\mn3;0}(n-1;n-2)$$
for all $1\leq j\leq n-2$, and
$h_{12\mn3;1}(n;n-1)=h_{12\mn3;1}(n-2)$, which means that
$$h_{12\mn3;1}(n)=h_{12\mn3;1}(n-1)+(n-1)h_{12\mn3;1}(n-2)+(n-2)h_{12\mn3;0}(n-3).$$

Case $r=2$: similarly as the above cases.
\end{proof}

\subsection{$\tau=13\mn2$}
\begin{theorem}\label{th32a}
  Let $r$ be a nonnegative integer. Then
$$H_{13\mn2;r}(x)=\frac{1-x}{1-2x}\delta_{r,0}+\frac{x^2}{1-2x}\sum\limits_{d=1}^r
\frac{H_{13\mn2;r-d}(x)-\sum\limits_{j=0}^{d-1}h_{13\mn2;r-d}(j)x^j}{(1-x)^d}.$$
\end{theorem}
\begin{proof}
Definitions imply $h_{13\mn2;r}(n;n)=h_{13\mn2;r}(n-1)$, and for
$1\leq j\leq n-1$,
$$h_{13\mn2;r}(n;j)=\sum\limits_{i=j}^{n-1} h_{13\mn2;r+i-j}(n-1;j).$$
By means of induction it is easy to obtain for $1\leq m\leq n-1$
$$h_{13\mn2;r}(n;n-m)=\sum\limits_{j=0}^{m-1}\binom{m-1}{j}h_{13\mn2;r-j}(n-1-m+j).$$
Now, if summing $h_{13\mn2;r}(n;n-m)$ over all $0\leq m\leq n-1$,
then we get
$$h_{13\mn2;r}(n)=h_{13\mn2;r}(n-1)+\sum\limits_{d=0}^r\sum\limits_{j=0}^{n-2-d}\binom{n-2-j}{d}h_{13\mn2;r-d}(j+d).$$
To find the desired result, we transfer the last equation to terms
of ordinary generating functions by use \cite[Lem.~7]{CM2}.
\end{proof}

In view of Theorem~\ref{th12} and Theorem~\ref{th32a} we have: the
number of permutations in $S_n(12\mn3)$ containing $13\mn2$
exactly $r$ times is the same number of permutations in
$S_n(21\mn3)$ containing $13\mn2$ exactly $r$ times. To verify
that by combinatorial bijective proof let $\pi$ be any
$12\mn3$-avoiding permutation; it is easy to see
$\pi=(\pi_1,\dots,\pi_{j-1},n,\pi')$ where
$\pi_1>\dots>\pi_{j-1}$, so the number of occurrences of $13\mn2$
in $\pi$ is given by $N:=n-1-\pi_{j-1}-(j-2)+N'$ where $N'$ the
number occurrences of $13\mn2$ in $\pi'$. On the other hand, let
$\beta=\Phi(\pi)$, so by definitions of $\Phi$ with induction
hypothesis (induction on length of $\pi$) we get that $\beta$
contains the same number $N$ of occurrences of $13\mn2$. Hence, by
means of induction we shall showed that $\Phi$ is a bijection, and
$H_{13\mn2;r}(x)=F_{13\mn2;r}(x)$ for all $r\geq 0$.

\subsection{$\tau=23\mn1$}
\begin{theorem}\label{th36}
Let $r$ be any nonnegative integer number. Then there exist a
polynomial $p_{r}(n)$ of degree at most $r$ with coefficient in
$\mathcal{Q}$ such that for all $n\geq r$
  $$h_{23\mn1;r}(n)=p_{r}(n)2^n.$$
\end{theorem}
\begin{proof}
Immediately, by definitions we have
\begin{lemma}\label{l361}
Let $n\geq 1$; then
$$\begin{array}{l}
h_{23\mn1;r}(n)=h_{23\mn1;r}(n-1)+\sum\limits_{j=1}^{n-1}h_{23\mn1;r}(n;j),\\
h_{23\mn1;r}(n;j)=\sum\limits_{i=j}^{n-1}h_{23\mn1;r-(j-1)}(n-1;i),\quad\mbox{for}\
1\leq j\leq n-1.
\end{array}$$
\end{lemma}

Hence, $h_{23\mn1;r}(n;n)=h_{23\mn1;r}(n;1)=h_{23\mn1;r}(n-1)$,
$h_{23\mn1;r}(n;j)=0$ for all $r+2\leq j\leq n-1$, and
$$h_{23\mn1;r}(n;j)=h_{23\mn1;r}(n;j-1)-h_{23\mn1;r+1-j}(n-1;j-1)$$
for all $2\leq j\leq r+1$.

Assume that $h_{23\mn1;d}(n)$ can be expressed as $p_d(n)2^n$ and
$h_{23\mn1;d}(n;j)$ can be expressed as $p_{d-1}(n)2^n$ where
$2\leq j\leq r+1$ for all $0\leq d\leq r-1$. The statement is
trivial for $r=0$, and by using the principle of induction with
the above explanation we get, immediately, the desired result.
\end{proof}

Theorem \ref{th36} with the initial values of the sequences
$h_{23\mn1;r}(n)$ for $r=0,1,2,3$ yield

\begin{corollary}
$$\begin{array}{l}
 {\rm (i)}\ \mbox{For all}\ n\geq1,\ \ h_{23\mn1;0}(n)=2^{n-1};\\
 {\rm (ii)}\ \mbox{For all}\ n\geq2,\ \ h_{23\mn1;1}(n)=(n-2)2^{n-3};\\
 {\rm (iii)}\ \mbox{For all}\ n\geq3,\ \
h_{23\mn1;2}(n)=(n-3)(n+8)2^{n-6};\\
{\rm (iv)}\ \mbox{For all}\ n\geq4,\ \
h_{23\mn1;3}(n)=1/3(n-4)(n^2+25n+42).
\end{array}$$
\end{corollary}

\subsection{$\tau=31\mn2$}
\begin{theorem}\label{th32b}
  Let $r$ be a nonnegative integer. Then
$$H_{31\mn2;r}(x)=\frac{1-x}{1-2x}\delta_{r,0}+\frac{x^2}{1-2x}\sum\limits_{d=1}^r
\frac{H_{31\mn2;r-d}(x)-\sum\limits_{j=0}^{d-1}h_{31\mn2;r-d}(j)x^j}{(1-x)^d}.$$
\end{theorem}
\begin{proof}
Definitions imply $h_{31\mn2;r}(n;1)=h_{13\mn2;r}(n-1)$, and for
$2\leq j\leq n-1$,
$$h_{31\mn2;r}(n;j)=\sum\limits_{i=j}^{n-1}h_{31\mn2;r}(n-1;i)=h_{31\mn2;r}(n-1)-\sum\limits_{i=1}^{j-1}h_{31\mn2;r}(n-1;i).$$
By means of induction it is easy to obtain for $1\leq m\leq n-1$
$$h_{31\mn2;r}(n;m)=\sum\limits_{j=0}^{m-1}(-1)^j\binom{m-1}{j}h_{31\mn2;r}(n-1-j).$$
Similarly as Theorem~\ref{th12} (or Theorem~\ref{th32a}), by using
the above equation with (it is easy to check by definitions)
$$h_{31\mn2;r}(n;n)=\sum\limits_{j=0}^r h_{31\mn2;r-j}(n-1;n-1-j),$$
we get the desired result.
\end{proof}

Again, we have $H_{31\mn2;r}(x)=F_{13\mn2;r}(x)$ for all $r\geq
0$. But here we failed to find a combinatorial explanation that
the number of permutations in $S_n(12\mn3)$ containing $13\mn2$
exactly $r$ times is the same number of permutations in
$S_n(21\mn3)$ containing $31\mn2$ exactly $r$ times.

\subsection{$\tau=32\mn1$}
\begin{theorem}\label{th37}
For all $n\geq 1$,
$$\begin{array}{l}
h_{32\mn1;0}(n)=\sum\limits_{j=0}^{n-2}(-1)^j\binom{n-1}{j+1}h_{32\mn1;r}(n-1-j)+\\
\qquad\qquad\qquad\qquad\qquad+\sum\limits_{j=1}^{r+1}\sum\limits_{i=0}^{j-1}(-1)^i\binom{j-1}{i}h_{32\mn1;r+1-j}(n-2-i).
\end{array}$$
\end{theorem}
\begin{proof}
Immediately, definitions yield
\begin{lemma}\label{l371}
Let $n\geq 1$; then
$$\begin{array}{l}
h_{32\mn1;r}(n)=h_{32\mn1;r}(n-1)+\sum\limits_{j=2}^{n}h_{32\mn1;r}(n;j),\\
h_{32\mn1;r}(n;j)=h_{32\mn1;r}(n-1)-\sum\limits_{i=1}^{j-1}h_{32\mn1;r}(n-1;i),
\quad\mbox{for}\ 1\leq j\leq n-1, \end{array}$$ and
    $$h_{32\mn1;r}(n;n)=h_{32\mn1;r}(n-1;1)+h_{32\mn1;r-1}(n-1;2)+\dots+h_{32\mn1;0}(n-1;r+1).$$
\end{lemma}

By means of induction with use Lemma~\ref{l371} we imply that for
all $1\leq m\leq n-1$
$$h_{32\mn1;r}(n;m)=\sum\limits_{i=0}^{m-1}(-1)^i\binom{m-1}{i}h_{32\mn1;r}(n-1-i).$$
On the other hand, using Lemma~\ref{l371} the third equality and
then using Lemma~\ref{l371} the first equality we get the desired
result.
\end{proof}

For example, Theorem \ref{th37} with \cite[Lem.~7]{CM2} (as
Example~\ref{exxx1}) yield the exact formula for $H_{32\mn1;r}(x)$
where $r=0,1$ (see Claesson and Mansour \cite{CM1} for the case
$r=0$).

\begin{corollary}
$$\begin{array}{l}
H_{32\mn1;0}(x)=\sum\limits_{k\geq0}\frac{x^{2k}}{p_{k-1}(x)p_{k+1}(x)};\\
H_{32\mn1;1}(x)=\sum\limits_{n\geq0}\left[\frac{x^2(1-(n+2)x)}{1-(n+1)x}\sum\limits_{k\geq0}\frac{u^{2(k+n)}}{p_{n+k}(x)p_{n+k+2}(x)}\right].
\end{array}$$
where $p_d(x)=\prod_{j=0}^d(1-dx)$.
\end{corollary}

\section{Further results}
The first possibility to extend the above result is to fix two
numbers of occurrences for two generalized patterns of type
$(2,1)$. For example, the number of permutations in $S_n$
containing $12\mn3$ exactly once and containing $13\mn2$ exactly
once is given by
$$(n^2-7n+14)2^{n-3}-2$$
for all $n\geq1$. Another example, the number of permutations in
$S_n$ containing $12\mn3$ exactly twice and containing $13\mn2$
twice is given by
$$(n^4-18n^3+163n^2-826n+1832)2^{n-7}-4n-14$$
for all $n\geq 1$. These results can be extended as follows.

\begin{theorem}\label{ext1}
Let us denote the number of permutations in $S_n$ containing
$12\mn3$ exactly $r$ times and containing $13\mn2$ exactly $s$
times by $a_n^{r,s}$; then there exists a polynomials $p(n)$ and
$q(n)$ of degree at most $r+s+1-\delta_{r,0}-\delta_{s,0}$ and
$r+s-\delta_{r,0}-\delta_{s,0}$, respectively, such that for all
$n\geq 1$
                    $$a_n^{r,s}=p(n)2^n+q(n).$$
\end{theorem}

Another direction to extend the results in above sections is to
restricted more than two patterns. For example, the number of
permutations in $S_n(12\mn3,13\mn2,21\mn3)$ is given by the
$(n+1)$th Fibonacci number (see \cite{CM1}). Again, this result
can be extended as follows.

\begin{theorem}\label{ext2}\ \\
(i) The ordinary generating function for the number of
permutations in $S_n(12\mn3$, $21\mn3)$ such containing $13\mn2$
exactly $r\geq1$ times is given by
$$\frac{x^2(1-x)^{r-1}}{(1-x-x^2)^{r+1}}$$
and for $r=0$ is given by $\frac{1}{1-x-x^2}$.\\
\ \\
(ii) The ordinary generating function for the number of
permutations in $S_n(12\mn3$, $21\mn3)$ such containing $23\mn1$
exactly $r\geq1$ times is given by
$$\frac{x^2(1-x)^{r-1}}{(1-x-x^2)^{r+1}}$$
and for $r=0$ is given by $\frac{1}{1-x-x^2}$.
\end{theorem}

In view of Theorem~\ref{ext2} suggests that there should exist a
bijection between the sets $\{12\mn3,21\mn3\}$-avoiding
permutations such containing $13\mn2$ exactly $r$ times and
$\{12\mn3,21\mn3\}$-avoiding permutations such containing $23\mn1$
exactly $r$ times for any $r\geq0$. However, we failed to produce
such a bijection, and finding it remains a challenging open
question.

\end{document}